\documentclass[24pt,twoside]{article}
\usepackage[centertags]{amsmath}
\usepackage{amsfonts}
\usepackage{amssymb}
\usepackage{amsthm}
\usepackage{newlfont}
\usepackage{color}
\usepackage[svgnames,dvipsnames]{xcolor}
\usepackage{graphicx}
\usepackage{pdfpages}

\date{\color{green}  2013 November}
\def \n {\noindent}
\usepackage{fancyhdr}
\usepackage{bm}
\usepackage[T1]{fontenc}
\usepackage[french]{babel}
\usepackage{array,multirow,makecell}
\setcellgapes{1pt}
\makegapedcells
\newcolumntype{R}[1]{>{\raggedleft\arraybackslash }b{#1}}
\newcolumntype{L}[1]{>{\raggedright\arraybackslash }b{#1}}
\newcolumntype{C}[1]{>{\centering\arraybackslash }b{#1}}

\begin{document}

\begin{center}
{\bf{\color{red}{\large  On the completeness of generalized eigenfunctions of the Hamiltonian of reggeon field theory in Bargmann space }}}\\
\end{center}

\begin{center}
{\bf Abdelkader Intissar$^{(1,2)}$}\\
\end{center}
\n {\it (1) Equipe d'Analyse spectrale, UMR-CNRS n: 6134, Université de Corse, Quartier Grossetti, 20 250 Corté-France.\\
E.mail: abdelkader.intissar@orange.fr\\}
\n {\it (2) Le Prador,129, rue du commandant Rolland, 13008 Marseille-France.\\}
\begin{center}
\n {\bf{\color{red}{\large  Abstract}}}
\end{center}

\n The aim of the present paper is to review and study some new spectral properties of the operator $H_{\mu, \lambda} = \mu A^{*}A + i\lambda A^{*}(A + A^{*})A$  which characterizes the Reggon field theory, where  $\mu$, $\lambda$ are real parameters and $i^{2} = -1$. $\displaystyle{A \longrightarrow \frac{d}{dz}}$ and $\displaystyle{A^{*} \longrightarrow z = x + i y ; (x, y) \in \mathbb{R}^{2}}$ respectively are boson annihilation and creation operators,  satisfying  the commutation relation $[A, A^{*}] = I$ in Bargmann space $\displaystyle{  \mathbb{B }= \{\varphi : \mathbb{C} \longrightarrow \mathbb{C}  \, entire \, ; \int_{\mathbb{C} }\mid \varphi(z)\mid^{2}e^{-\mid z \mid^{2}}dx dy < \infty \}}$ with standard scalar product $\displaystyle{ < \varphi, \psi > = \frac{1}{\pi} \int_{\mathbb{C} } \varphi(z) \overline{\psi(z)}e^{-\mid z \mid^{2}}dx dy }$\\

\n The Hamiltonian takes the form $\displaystyle{H_{\mu, \lambda} = \mu z \frac{d}{dz} + i\lambda(z \frac{d^{2}}{dz^{2}} + z^{2} \frac{d}{dz}): = \mu \mathcal{N} +  i\lambda H_{I}}$,  where $\displaystyle{\mathcal{N} = z \frac{d}{dz}}$ and  $\displaystyle{ H_{I} = z \frac{d^{2}}{dz^{2}} + z^{2} \frac{d}{dz}}$.\\

\n The Hamiltonian $H_{\mu,\lambda}$ is non-Hermitian with respect to the above standard scalar product. Hence the question arises, whether the eigenfunctions by the finite norm condition  form a complete basis. \\

\n The main new results of this article are the determination of the boundary conditions for the eigenvalue problem associated to  $H_{\mu,\lambda}$  and the proof of the completeness of the basis mentioned in the question above by transforming $H_{\mu,\lambda}$ in the Hermitian form.\\

\n This article also provides an explicit integral form of the inverse of $H_{\mu,\lambda}$ on the negative imaginary axis $z = -iy;  \,  y > 0$ and some spectral properties of $H_{\mu,\lambda}^{-1};  \, \mu > 0$.\\

\n {\bf Keywords : } Spectral theory ; Gribov operators ; semigroups ; Non-self-adjoint operators ; Bargmann space ; Reggeon field theory.\\

\n $\underline{\quad\quad\quad\quad\quad\quad\quad\quad\quad\quad\quad\quad\quad\quad\quad\quad\quad\quad \quad\quad\quad}$\\

\n {\bf{\color{red}  $\S$\,1. Introduction}}\\

\n Le $\mathbb{B}$ be the Bargmann space {\color{blue}[4]} defined by :\\

\n $\displaystyle{  \mathbb{B} = \{\varphi : \mathbb{C} \longrightarrow \mathbb{C}  \, entire \, ; \int_{\mathbb{C} }\mid \varphi(z)\mid^{2}e^{-\mid z \mid^{2}}dx dy < \infty \}}$ \hfill { } {\color{blue} (1.1)} \\ 

\n with standard scalar product \\

\n $\displaystyle{ < \varphi, \psi > = \frac{1}{\pi}\int_{\mathbb{C} } \varphi(z) \overline{\psi(z)}e^{-\mid z \mid^{2}}dx dy }$  \hfill { } {\color{blue} (1.2)} \\ 

\n and  with usual basis\\

\n $\displaystyle{e_{n}(z) = \frac{z^{n}}{\sqrt{n!}}; \,  n \in \mathbb{N}}$ \hfill { } {\color{blue} (1.3)}\\

\n and\\

\n $\mathbb{P}$  is the set of polynomials which is dense in $\mathbb{B}$.\\

\n The constant $\displaystyle{\frac{1}{\pi}}$ is chosen in scalar product such that  the norm of the constant function $\varphi(z) = 1$ is one. \\

\n Let $\displaystyle{ \mathbb{B}_{s} = \{(a_{n})_{n=0}^{\infty} \in \mathbb{C} ; \sum_{n=0}^{\infty}n!\mid a_{n}\mid^{2} < \infty \}}$ $\hfill { } {\color{blue}(1.4)}$\\ 

\n $\mathbb{B}_{s}$ is also a Hilbert space. The scalar product on $\mathbb{B}_{s}$ is defined by\\
 
\n  $\displaystyle{ < (a_{n}) , (b_{n} >_{s} = \sum_{n=0}^{\infty}n!a_{n}\overline{b_{n}}}$ $\hfill { } {\color{blue}(1.5)}$\\
  
\n  and the associated norm is denoted by $\mid\mid . \mid\mid_{s}$.\\

\n It is well known that  $\mathbb{B}$ is related to $\mathbb{B}_{s} $ by unitary transform of  $\mathbb{B}$ onto $\mathbb{B}_{s} $, given by the following transform :\\

\n $\displaystyle{ I: \mathbb{B} \longrightarrow \mathbb{B}_{s} , \varphi(z) = \sum_{n=0}^{\infty}a_{n}z^{n} \longrightarrow I(\varphi) = (\frac{1}{n!}\varphi^{(n)}(0))_{n=0}^{\infty} = (a_{n})_{n=0}^{\infty}}$ $\hfill { } {\color{blue}(1.6)}$\\
                             
\n In Bargmann's representation, we have : \\

\n the annihilation operator $ A$ is :\\

\n $\displaystyle{A = \frac{d}{dz}}$ with domain  $\displaystyle{D(A) = \{\varphi \in \mathbb{B}; \frac{d}{dz} \varphi \in \mathbb{B}\}}$ \hfill { } {\color{blue} (1.7)}\\

\n  the creation operator $A^{*}$ is : \\

\n the multiplication by $z$ with domain  $\displaystyle{D(A^{*}) = \{\varphi \in \mathbb{B}; z\varphi \in \mathbb{B}\}}$  \hfill { } {\color{blue} (1.8)}\\

\n It is well known that  $A^{*}$ is adjoint of $A$ (with respect above standard scalar product), $D(A)$ is dense in $\mathbb{B}$, $D(A) = D(A^{*})$ and the injection  $D(A)$ $\longrightarrow$ $\mathbb{B}$ is a compact  mapping.(see  {\color{blue}[13]} or lemma 0.2 in {\color{blue}[15]} for an elementary proof).\\

\n and the Hamiltonian of reggeon field theory in zero transverse dimensions takes the following  form :\\

\n {\color{red} ($\bullet_{1} $)} $\displaystyle{H_{\mu, \lambda} = \mu z \frac{d}{dz} + i\lambda(z \frac{d^{2}}{dz^{2}} + z^{2} \frac{d}{dz}): = \mu \mathcal{N} +  i\lambda H_{I}}$ ; $(\mu, \lambda) \in \mathbb{R}^{2}$ and $i^{2} = - 1$ \hfill { } {\color{blue} (1.9)} \\

\n with maximal domain $\displaystyle{D(H_{\mu, \lambda}) = \{\varphi \in \mathbb{B}; H_{\mu, \lambda} \varphi \in \mathbb{B}\}}$ .\\

\n where \\

\n $\displaystyle{\mathcal{N} = z \frac{d}{dz}}$  with domain $\displaystyle{D(\mathcal{N}) = \{\varphi \in \mathbb{B}; \mathcal{N} \varphi \in \mathbb{B}\}}$ \hfill { } {\color{blue} (1.10)} \\

\n and \\

\n  $\displaystyle{ H_{I} = z \frac{d^{2}}{dz^{2}} + z^{2} \frac{d}{dz}}$  with maximal domain $\displaystyle{D(H_{I}) = \{\varphi \in \mathbb{B}; H_{I} \varphi \in \mathbb{B}\}}$.  \hfill { } {\color{blue} (1.11)} \\

\n {\color{red} ($\bullet_{2} $)} If we call generalized Laguerre-type derivatives the operators of the form $\displaystyle{\mathcal{D}(r, s) = A^{r}\mathcal{N}^{s}}$, with $r = 1, 2, ... , s = 0, 1, ... $ then as $\mathcal{N}^{s}$  conserves the number of bosons, we observe that these operators act as monomials in boson operators which annihilate $r$ bosons and that the operator $H_{\mu, \lambda}$ can be written  in the form\\

\n $H_{\mu, \lambda} = \mu A^{*}A + i\lambda (\mathcal{D}(1, 1) + \mathcal{D}^{*}(1, 1) - (A + A^{*}))$ \hfill { } {\color{blue} (1.12)} \\

\n {\color{red} ($\bullet_{3} $)} In 1987, we have given  in {\color{blue}[13]} many spectral properties of $H_{\mu, \lambda}$ for $\mu > 0$ in particular the minimal domain of $\displaystyle{H_{\mu, \lambda}}$ coincides with its maximal domain, the positiveness of its eigenvalues , the existence of the smallest eigenvalue $\sigma_{0} \neq 0$ and an asymptotic expansion of its semigroup $\displaystyle{e^{-tH_{\lambda, \mu}} }$ as $t \longrightarrow + \infty$.\\

\n Where the minimal domain of $\displaystyle{H_{\mu, \lambda}}$ is:  \\

\n $\displaystyle{ D(H_{\mu, \lambda}^{min}) = \{ \varphi \in \mathbb{B} , \exists \, p_{n} \in \mathbb{P} , p_{n} \longrightarrow \varphi , \exists \, \psi \in \mathbb{B} ; H_{\mu, \lambda}p_{n} \longrightarrow \psi \}}$ .\hfill { } {\color{blue} (1.13)} \\\

\n {\color{red} ($\bullet_{4} $)} In 1984, Ando and Zerner have schown in {\color{blue}[3]} that the smallest eigenvalue $\sigma_{0}$ of $H_{\mu, \lambda}$ can be  analytically continued  in $\mu$ on entire real axis.\\

\n {\color{red} ($\bullet_{5} $)} In 1998, we have given  in {\color{blue}[14]} the boundary conditions at infinity for a description of all maximal dissipative extensions in Bargmann space of the minimal Heun's operator $\displaystyle{ H_{I} = z \frac{d^{2}}{dz^{2}} + z^{2} \frac{d}{dz}}$. The characteristic functions of the dissipative extensions have computed and some completeness theorems have obtained for the system of generalized eigenvectors. it is well known that the restriction $H_{I}^{min}$ of the closure of $H_{I}$ on the polynomials set $\mathbb{P}$ is symmetric.\\

\n  But  the minimal domain $\displaystyle{ D(H_{I}^{min}) = \{ \varphi \in \mathbb{B} , \exists \, p_{n} \in \mathbb{P} , p_{n} \longrightarrow \varphi , \exists \, \psi \in \mathbb{B} ; H_{I}p_{n} \longrightarrow \psi \}}$ of $\displaystyle{ H_{I}}$ is different of its maximal domain $D(H_{I}) = \{ \varphi \in \mathbb{B} ; H_{I}\varphi \in \mathbb{B} \}$. \\

\n {\color{red} ($\bullet_{6} $)} It is also well known that $\displaystyle{ H_{I}}$ is chaotic operator in Devaney's sense see {\color{blue}[8]}  or {\color{blue}[16]}, in particular its spectrum is $\sigma(H_{I}) = \mathbb{C}$.\\

\n {\color{red} ($\bullet_{7} $)} The Hamiltonian $H_{\mu,\lambda}$ is non-Hermitian with respect to above standard scalar product and the domain of the adjoint and anti-adjoint parts are not included in each other, nor is the domain of their commutator.\\

\n {\color{red} ($\bullet_{8} $)} Hence the question arises, whether the eigenfunctions by the finite norm condition  form the complete basis.  The aim new results of this article is the determination of boundary conditions of eigenvalue problem associated to  $H_{\mu,\lambda}$  and to give a answer to above question by transform $H_{\mu,\lambda}$  to the Hermitian form. Also, in this article we give an explicit integral form of the inverse of $H_{\mu,\lambda}$ on the negative imaginary axis $z = -iy, \, y >0$  and some spectral properties of $H_{\mu,\lambda}^{-1}; \mu > 0$ on a weighted space.\\

\n Now we give an outline of the content of this paper:\\

\n  The Section 2 deals the action of some elementary operators on Bargmann space usefull to study of eigenvalue problem associated to the Hamiltonian $H_{\mu, \lambda}$.\\

\n In section 3, we give an asymptotic analysis of the solutions of bi-confluent Heun equation associated to eigenvalue problem of the Hamiltonian $H_{\mu, \lambda}$. \\

\n In section 4,  we give on the negative imaginary axis $z = -iy , y > 0$ a transformation procedure of $H_{\mu, \lambda}$ to a symmetric operator with compact resolvent and we show the main result of this paper that the eigenfunctions of $H_{\mu, \lambda}$ form  a complete basis.\\

\n And\\

\n In section 5, we give an explicit integral form of the inverse of $H_{\mu,\lambda}$ on the negative imaginary axis $z = -iy;  \,  y > 0$ and some spectral properties of $H_{\mu,\lambda}^{-1};  \, \mu > 0$ as integral operator on a weighted Hilbert space.\\
 
\n {\bf{\color{red}$\S$  2 \,  Action of some elementary operaors on Bargmann space usefull to study of eigenvalue problem associated to the Hamiltonian $H_{\mu, \lambda}$}}\\

\n {\color{red}$\bullet$} We introduce the operator $P$  (changing the sign of $z$) defined by\\

\n $\displaystyle{PzP^{-1} = -z}$ and $\displaystyle{P\frac{d}{dz}P^{-1} = - \frac{d}{dz}}$  then we deduce that\\

\n  $\displaystyle{PHP^{-1} =  \mu z \frac{d}{dz} - i\lambda (z\frac{d^{2}}{dz^{2}} + z^{2} \frac{d}{dz}): = H_{\mu, -\lambda}}$ (the formal adjoint of $H_{\mu, \lambda}$) \hfill { } {\color{blue} (2.1)} \\

\n The operator $P$ is unitary and Hermitian in Bargmann space, since \\

\n $\displaystyle{ < P\varphi , P \psi > =   < \varphi ,  \psi > }$ and $\displaystyle{ <\varphi , P \psi >  =  < P\varphi , \psi > }$\\

\n Now, let the following  eigenvalue problem :\\

\n $H_{\mu, \lambda}\varphi_{n} = \sigma_{n} \varphi_{n}$ with real $\sigma_{n}$ \hfill { } {\color{blue} (2.2)} \\

\n and \\

\n $\displaystyle{\varphi_{n}^{*}(z) = - P\varphi_{n}(z) = - \varphi_{n}(-z)}$ \hfill { } {\color{blue} (2.3)} \\

\n  then  we get:\\

\n $\displaystyle{ <\varphi_{m}^{*}, H_{\mu, \lambda}\varphi_{n} > = \sigma_{n} <\varphi_{n}^{*},  \varphi_{n} >}$  \hfill { } {\color{blue} (2.4)} \\

\n and\\

\n $\displaystyle{ H_{\mu, -\lambda}\varphi_{m}^{*} = - H_{\mu, -\lambda}P \varphi_{m} = - PH_{\mu, \lambda}\varphi_{m} =\sigma_{m} \varphi_{m}^{*}}$ \hfill { } {\color{blue} (2.5)} \\

\n and \\

\n  $\displaystyle{ < H_{\mu, -\lambda} \varphi_{m}^{*} ; \varphi_{n} > =  \sigma_{m} < \varphi_{m}^{*} ; \varphi_{n} > }$  \hfill { } {\color{blue} (2.6)} \\

\n it follows that \\

\n $\displaystyle{(\sigma_{n} - \sigma_{m})< \varphi_{m}^{*}, \varphi_{n} > = < \varphi_{m}^{*}, H_{\mu, \lambda}\varphi_{n} > - <  H_{\mu, -\lambda}\varphi_{m}^{*} 
, \varphi_{n} > }$\\

\n $\displaystyle{ = < \varphi_{m}^{*}, H_{\mu, \lambda}\varphi_{n} >  - < \varphi_{m}^{*}, H_{\mu, \lambda}\varphi_{n} >}$ $\displaystyle{  =  0}$\\

\n so that \\

\n $\displaystyle{< \varphi_{m}^{*}, \varphi_{n} > = 0}$ when $\sigma_{n} \neq \sigma_{m}$ \hfill { } {\color{blue} (2.7)} \\

\n {\color{red}$\bullet$} Determination of an orthogonal scalar product in Bargmann space associated to the operator $P$.\\

\n We introduce the operator $\nu$ (the Hermitian operator of the sign)  on the eigenfunctions $\{\varphi_{n}\}, n = 1, 2, ..;$  of $H_{\mu, \lambda}$ defined by\\

\n $\displaystyle{\nu( \varphi_{n} ) = \nu_{n}\varphi_{n}}$\\

\n where\\

\n $\left \{ \begin{array} {c} \nu_{n} = + 1, \quad if  < \varphi_{n} , P\varphi_{n} > {\color{blue}>  0} \\
\quad\\
\nu_{n} = -1,\quad  if  < \varphi_{n} , P\varphi_{n} >  {\color{blue} < 0} \\
\end{array} \right.$\\
\quad \\

\n We introduce  the  positively defined scalar product  $\displaystyle{<< \varphi , \psi >>  =  < \varphi , \nu P\psi > }$ with respect to which the eigenstates of $H_{\mu, \lambda}$ are orthogonal.\\

\n {\bf {\color{red} $\S$ \, 3 Asymptotic analysis of the solutions of bi-confluent Heun equation associated to eigenvalue problem of the Hamiltonian $H_{\mu, \lambda}$}}\\

\n We recall that  $H_{\mu, \lambda}$ has the form $\displaystyle{p(z)\frac{\partial^{2}}{\partial z^{2}} + q(z)\frac{\partial}{\partial z} }$ with $p(z) = i\lambda z$ and $q(z) =  i\lambda z^{2} + \mu z$ of degree {\color{red}$2$} and $H_{\mu, \lambda}$ is of Heun operator type.\\

\n The eigenvalue problem associated to $H_{\mu, \lambda}$ does not satisfy the classical ordinary differential equations of the form:\\

\n $\displaystyle{p(z)\frac{\partial^{2}\varphi}{\partial z^{2}} + q(z))\frac{\partial \varphi}{\partial z}  = \sigma \varphi}$ $\hfill { } {\color{blue}(3.1)}$\\

\n where $p(z)$ is a polynomial of degree  at most two, $q(z)$ is a polynomial of the degree exactly {\color{red}one} and $\sigma$ is a constant.\\

\n {\color{red}$\bullet$} If $\lambda \neq  0$, the eigenvalue problem associated to $H_{\mu, \lambda}$ can be written as follows :\\

\n $\displaystyle{z\frac{\partial^{2}\varphi}{\partial z^{2}} +( z^{2} - i\rho z)\frac{\partial \varphi}{\partial z}  = - i\frac{\sigma}{\lambda} \varphi; \quad \rho = \frac{\mu}{\lambda}}$ $\hfill { } {\color{blue}(3.2)}$\\

\n Let $ z = i\sqrt{2}\xi$ and $\varphi(\sqrt{2}\xi) = \psi(\xi)$ then {\color{blue}(3.2)} can be transformed to :\\

\n $\displaystyle{\xi\frac{\partial^{2}\psi}{\partial \xi^{2}} +(- 2\xi^{2}  + \rho\sqrt{2} \xi)\frac{\partial \psi}{\partial\xi} - \frac{\sigma\sqrt{2}}{\lambda} \psi  = 0}$  $\hfill { } {\color{blue}(3.3)}$\\

\n the above equation belongs to bi-confluent Heun equations, denoted by $BHE(\alpha, \beta, \gamma, \delta)$ which are in the form :\\

\n $\displaystyle{ \xi \frac{d^{2}u}{d\xi^{2}}  + ( 1 + \alpha  - \beta \xi - 2\xi^{2}) \frac{du}{d\xi} + \{(\gamma - \alpha - 2)\xi - \frac{1}{2}(\delta + (\alpha+1)\beta)\}u = 0}$; \, $(\alpha, \beta, \gamma, \delta) \in \mathbb{C}^{4}$  $\hfill { } {\color{blue}(3.4)}$\\

\n It is well known that \\

\n {\color{red}$\bullet$} the singular points of the bi-confluent Heun equation {\color{blue}(3.4)} are $0$ and $\infty$. \\

\n {\color{red}$\bullet$} The singularity is regular at $0$ and irregular at  $\infty$.\\

\n {\color{red}$\bullet$} At $0$, if $\displaystyle{(\alpha, \beta, \gamma, \delta) \in (\mathbb{C} - \mathbb{Z})\times \mathbb{C}^{3}}$ then a basis of $BHE(\alpha, \beta, \gamma, \delta)$ is given by \\

\n $\displaystyle{\{ N(\alpha, \beta, \gamma, \delta, \xi), \xi^{-\alpha}N(-\alpha, \beta, \gamma, \delta, \xi)\}}$ where  $\displaystyle{N(\alpha, \beta, \gamma, \delta, \xi) = \Gamma(\alpha)\sum_{n= 0}^{\infty}\frac{A_{n}(\alpha, \beta, \gamma, \delta)}{\Gamma(\alpha + 1 + n)n!}\xi^{n}}$, \\
 
 \n the $A_{n}$ are polynomials in $\alpha, \beta, \gamma, \delta$ defined by the following relation\\

\n $\displaystyle{A_{n+2} = [\beta(n + 1) + \frac{1}{2}(\delta + \beta(1 + \alpha))]A_{n+1} - (\gamma - 2 - \alpha - 2n)(n+ 1)(n + 1 + \alpha)A_{n}}$ for  $ n \geq 0$ and $A_{-1} = 0, A_{0} = 1$, (see proposition 5 of ref. {\color{blue}[22]} or ref.  {\color{blue}[5]}) \hfill { } {\color{blue} ($\star_{1}$)}  \\

\n {\color{red}$\bullet$} If $ \alpha \in \mathbb{Z}$ there are logarithmic terms in the bases of solutions.of $BHE(\alpha, \beta, \gamma, \delta)$ at $0$ which is a regular singularity (except for some values see {\color{blue}[20]}).  \hfill { } {\color{blue} ($\star_{2}$)}  \\

\n {\color{red}$\bullet$}  If $\displaystyle{(\alpha, \beta, \gamma, \delta) \in \mathbb{C}^{4}}$ then $BHE(\alpha, \beta, \gamma, \delta)$ admits as basis of formal solutions at $\infty$,\\

\n  $\displaystyle{\{\xi^{-\frac{1}{2}(\gamma - \alpha - 2)} \sum_{n=0}^{\infty} a_{n}(\alpha, \beta, \gamma, \delta)\xi^{-n}, e^{-\frac{1}{2} (\gamma + \alpha + 2} e^{\xi^{2} + \beta \xi} \sum_{n=0}^{\infty} a_{n}(\alpha, i\beta, -\gamma, -i\delta)(i\xi)^{-n}\}}$.\\

\n where the complex numbers an are defined by \\

\n  $\displaystyle{2(n+2)a_{n+2} = [\frac{1}{2}(\delta + \beta(\gamma -1) - \beta(n + 1)]a_{n + 1} - (\frac{\gamma - \alpha - 2}{2} - n )( \frac{\gamma + \alpha - 2}{2} - n)a_{n}}$  for  $ n \geq 0$ and $a_{-1} = 0, a_{0} = 1$.  (see proposition $8$ of ref  {\color{blue}[22]} ) \hfill { } {\color{blue} ($\star_{3}$)}  \\

\n These equations are more complicated than the equations for classical special functions such as Bessel functions, hypergeometric functions, and confluent hypergeometric functions (see {\color{blue}[7]} and {\color{blue}[1]}).\\

\n {\bf {\color{red} Remark 3.1}}\\

\n ( i) If we take  $\displaystyle{\alpha = -1, \beta = -\rho \sqrt{2}, \gamma = 1}$ and $\displaystyle{\delta = \frac{2\sigma\sqrt{2}}{\lambda}}$ in {\color{blue}(3.4)} then we get {\color{blue}(3.3)}, so it is \\ $\displaystyle{BHE(-1, -\frac{\mu \sqrt{2}}{\lambda}, 1, -\frac{2\sigma\sqrt{2}}{\lambda})}$ .\\

\n (ii) As $\alpha = -1 \in \mathbb{Z} $ then the proposition $5$ of {\color{blue}[22]} (see {\color{blue} $\star_{1}$}) is not applicable.\\

\n (iii) We will show there are logarithmic terms in the bases of solutions of $\displaystyle{BHE(-1, -\frac{\mu \sqrt{2}}{\lambda}, 1, \frac{2\sigma\sqrt{2}}{\lambda})}$ at $0$  which is a regular singularity  by using the classical Frobenius-Fuch's method. \hfill { } {\color{blue} $\square$}\\

\n In {\color{blue}(3.2)} denoting $\displaystyle{  \frac{d}{dz}\varphi (z) = \varphi^{'}(z)}$  and $\displaystyle{ \frac{d^{2}}{dz^{2}}\varphi (z) = \varphi^{"}(z) }$ where a prime throughout denotes a derivative with respect to $z$, then the eigenvalue problem associated to $H_{\mu, \lambda}; \lambda \neq 0$ can be written as follows \\

\n $\displaystyle{\varphi^{"}(z) + p(z) \varphi^{'}(z) + q(z) \varphi(z) = 0};\quad  \displaystyle{ p(z) = (z - i\rho)}$ and $\displaystyle{ q(z) =  \frac{i\sigma}{\lambda z}}$ $\hfill { } {\color{blue}(3.5)}$\\ 

\n Eq. {\color{blue} (3.5)} is bi-confluent equation designated  in  the  Ince classification {\color{blue}[11]} as {\bf{\color{red}[$0, 1, 1_{4}$]}} which has, as  in   the   case   of   the   confluent  hypergeometric equation, two singularities located, respectively, at the origin and infinity: one regular,  and  the  second  irregular  with  a  singularity  rank  higher  by  unity  than  that  for  the  confluent hypergeometric equation.\\

\n We recall that $z_{0}$ is regular singular point of {\color{blue}(3.5)} if\\

\n {\bf $\bullet$} $z_{0}$ is singular;\\

\n {\bf $\bullet$} $\displaystyle{(z- z_{0}) p(z) }$ and $(z- z_{0})^{2} q(z)$ are analytic at $z_{0}$. $\hfill { } {\color{blue}(3.6)}$\\

\n As $\displaystyle{lim \, zp(z) = lim \, z^{2}q(z) = 0}$ as $z \longrightarrow 0$ then  $z_{0} = 0$ is regular singular point of {\color{blue}(3.5)}.\\

\n The usual power series method, that is setting $\displaystyle{\varphi(z) = \sum_{n=0}^{+\infty}a_{n}(z - z_{0})^{n}}$, breaks down if $z_{0}$ is a singular point. Here "breaks down" means "cannot find all solutions".\\

\n The  solution  of  the  our bi-confluent  Heun  equation  in  a  power  series  of  $z$  can  be  easily  constructed by applying the Method of Frobenius to {\color{blue}(3.5)} whose the basic idea is to look for solutions of the form $\displaystyle{\varphi(z) = z^{s}\sum_{n=0}^{+\infty}a_{n}z^{n}}$ :\\

\n Substitute the expansion $\displaystyle{\varphi(z) = z^{s}\sum_{n=0}^{+\infty}a_{n}z^{n}}$ into the equation {\color{blue}(3.5)}, we have\\

\n $(\displaystyle{ z^{s}\sum_{n=0}^{+\infty}a_{n}z^{n}})^{"}$ + $(z - i\rho) (\displaystyle{ z^{s}\sum_{n=0}^{+\infty}a_{n}z^{n}})^{'}$ +  $\displaystyle{\frac{i\sigma}{\lambda z}}$ $\displaystyle{z^{s}\sum_{n=0}^{+\infty}a_{n}z^{n+s} = 0}$ $\hfill { } {\color{blue}(3.7)}$\\

\n Carrying out the differentiation, we reach\\

\n  $\displaystyle{\sum_{n=0}^{+\infty}(n+s)(n+s -1)a_{n}z^{n+s-2} \quad + \quad  \sum_{n=0}^{+\infty}(n+s)a_{n}z^{n+s} \quad +\quad}$ \\
$\displaystyle{  i(\frac{\sigma}{\lambda} - \rho)\sum_{n=0}^{+\infty}(n+s)a_{n}z^{n+s-1} = 0} \hfill { } {\color{blue}(3.8)}$\\

\n Shifting index:\\

\n {\bf $\bullet$} $\displaystyle{\sum_{n=0}^{+\infty}(n+s)a_{n}z^{n+s-1} = \sum_{n=1}^{+\infty}(n+s-1)a_{n-1}z^{n+s-2}}$$\hfill { } {\color{blue}(3.9)}$\\

\n {\bf $\bullet$} $\displaystyle{\sum_{n=0}^{+\infty}(n+s)a_{n}z^{n+s} = \sum_{n=1}^{+\infty}(n+s-3)a_{n-2}z^{n+s-2}; a_{-1} = 0}$$\hfill { } {\color{blue}(3.10)}$\\

\n Now if we put $\alpha = i(\frac{\sigma}{\lambda} - \rho)$ then the equation becomes\\

\n $\displaystyle{s(s-1)a_{0}z^{s-2} \quad +\quad \sum_{n=1}^{+\infty}[(n+s)(n+s -1)a_{n} \quad + \quad \alpha (n+s-1)a_{n-1} \quad + }$ \\

\n $\displaystyle{ (n+s-3)a_{n-2}]z^{n+s-2} = 0}$ $\hfill { } {\color{blue}(3.11)}$\\

\n and the indicial equation is \\

\n $s(s-1) = 0$.$\hfill { } {\color{blue}(3.12)}$\\

\n by using the theorem of Fuchs : \\

\n {\bf {\color{red}Theorem 3.2}} (theorem of  Fuchs)\\

\n Consider  the equation  $\displaystyle{\varphi"(z) + p(z)\varphi^{'} + q(z)\varphi(z) = 0}$\\

\n Let $z_{0}$ be a regular singular point. That is $\displaystyle{(z- z_{0})p(z) = \sum_{n=0}^{+\infty}p_{n}(z- z_{0})^{n}}$ and \\

\n  $\displaystyle{(z- z_{0})^{2}q(z) = \sum_{n=0}^{+\infty}q_{n}(z- z_{0})^{n}}$ with certain radii of convergence.\\

\n Let $r$ be no bigger than the radius of convergence of either $(z- z_{0})p(z)$ or $(z- z_{0})^{2}q(z)$ and $s_{1}, s_{2}$ solve the indicial equation \\

\n $\displaystyle{s (s -1) + p_{0} s + q_{0} =0}$.$\hfill { } {\color{blue}(3.13)}$\\

\n Then\\

\n 1. If $s_{1} \neq s_{2}$  and $s_{1} - s_{2}$ is not an integer, then the two linearly independent solutions are given by\\

\n $\displaystyle{\varphi_{1}(z) = (z- z_{0})^{s_{1}}\sum_{n=0}^{+\infty}a_{n}(z- z_{0})^{n}}$,  $\displaystyle{\varphi_{2}(z) = (z- z_{0})^{s_{1}}\sum_{n=0}^{+\infty}b_{n}(z- z_{0})^{n}}$$\hfill { } {\color{blue}(3.14)}$

\n 2. If $s_{1} = s_{2}$ then $\varphi_{1}(z)$ is given by the same formula as above, and $ \varphi_{2}(z)$ is of the form\\

\n $\displaystyle{\varphi_{2}(z) =\varphi_{1}(z) ln(z-z_{0}) + (z- z_{0})^{s_{1}} \sum_{n=0}^{+\infty}d_{n}(z- z_{0})^{n}}$$\hfill { } {\color{blue}(3.15)}$\\

\n 3. If  $s_{1} - s_{2}$ is an integer, then take $s_{1}$ to be the larger root (More precisely, when $s_{1}$, $s_{2}$ are both complex, take $s_{1}$ to be the one with larger real part, that is $\Re e s_{1} > \Re e s_{2}$). Then $\varphi_{1}(z) $ is still the same, while\\

\n $\displaystyle{\varphi_{2}(z) = c\varphi_{1}(z) ln(z-z_{0}) + (z- z_{0})^{s_{2}} \sum_{n=0}^{+\infty}e_{n}(z- z_{0})^{n}}$ $\hfill { } {\color{blue}(3.16)}$ \\

\n (Note that c may be $0$).\\

\n All the solutions constructed above converge at least for $0 < \mid z - z_{0}\mid < r$. \hfill { } {\color{blue} $\square$}\\

\n The proof of this theorem is through careful estimate of the size of an using the recurrence relation, see {\color{blue}[2]} :\\

\n In fact the converse of this theorem is also true. That is if all solutions of the equation satisfies $\displaystyle{Lim (z - z_{0})^{s}\varphi(z) = 0}$ as $z \longrightarrow z_{0}$ for some $s$, then $(z - z_{0})p(z)$ and $(z - z_{0})^{2}q(z)$ are analytic at $z_{0}$. This is called FuchsTheorem. Its proof is a tour de force of complex analysis and can be found in {\color{blue}[23]} see also, the following references about this topic {\color{blue}[10]}  or  {\color{blue}[21]}  \\

\n  We deduce that {\color{blue}(3.5)} has tow  linearly independent solutions given by\\

\n $ \left \{ \begin{array} [c] {l}\n  \displaystyle{\varphi_{1}(z) = z \eta(z) ; \quad  \eta \quad is \quad  analytic , \eta'(0) \neq 0}\\
\quad\\
\n \displaystyle{\varphi_{2}(z) =  c zln(z)\xi(z)  +  \sum_{n=0}^{+\infty}c_{n}z^{n}; \quad \xi(z) \quad is \quad analytic , \xi'(0) \neq 0}\\
\end{array} \right .$$\hfill { } {\color{blue}(3.17)}$
\quad\\

\n As $\varphi_{2}(z)$  contains non-analytic terms then it is not acceptable  in Bargmann space and thus it is rejected then  \\

\n $\varphi (z) \thicksim z$ \quad  ($ z  \longrightarrow 0$) $\hfill { } {\color{blue}(3.18)}$\\

\n Now, we turn to the point $z = \infty$, which is the second possible singular point for the solution of {\color{blue}(3.5)}). To analyse the point $z = \infty$, we can first perform the change of independent variable from $z$ to $\xi$, $\displaystyle{z = \frac{1}{\xi}}$ and $\varphi(z) = \varphi(\frac{1}{\xi}) = \phi(\xi)$ and study the behaviour of the transformed equation at $\xi = 0$ to get :\\

\n $\left \{ \begin{array} [c] {l}\displaystyle{ \frac{d}{dz} = -\xi^{2}\frac{d}{d\xi}\,   and\,  \frac{d^{2}}{dz^{2}} = \xi^{4}\frac{d^{2}}{d\xi^{2}} + 2\xi^{3}\frac{d}{d\xi}}\\
\quad\\
\n \displaystyle{\phi^{"}(\xi) +  p(\xi) \phi^{'}(\xi)  + q(\xi)\phi(\xi) = 0} \\
\quad\\
\n where\\
\quad\\
\n \displaystyle{p(\xi) = - \frac{2\xi^{2} - i\rho \xi + 1}{\xi^{3}}}\, and\,  \displaystyle{q(\xi) = \frac{i\sigma}{\lambda \xi^{3}}}\\
\end{array} \right .$$\hfill { } {\color{blue}(3.19)}$\\
\quad\\

\n  {\color{red}$\bullet$} The origin  is regular singular point of {\color{blue}(3.19)} if  $\displaystyle{\xi p(\xi)}$  and  $\displaystyle{ \xi^{2}q(\xi)}$ are analytic functions in the neighborhood of the origin. As these latter conditions are not fulfilled, then the origin is an irregular singular point and it is of order $3$.\\

\n  Consequently, the theorem of Fauchs is not applicable to {\color{blue}(3.19)} and this equation cannot be approximated by an Euler equation in the neighborhood of the origin.\\ 

\n If we like to explore some of these techniques, a good starting point would be Chapter $3$ of {\color{blue}[6]} or the book of Kristensson {\color{blue}[19]}.\\

\n The case of an irregular singular point is much more difficult to address, although some techniques do exist for obtaining the solution to a homogeneous second-order linear differential equation in the vicinity of an irregular singular point.\\

\n The applications of the concepts on irregular singular points in the current literature to our second order linear differential equation does not inform us about the solutions belonging to the domains of our operators, in particular to the Bargmann space and what are the boundary conditions verified by them at infinity ?\\
 \n For example, if we look for a solution of {\color{blue}(3.19)} in the form $\displaystyle{\phi(\xi) = \sum_{k=0}^{+\infty}a_{k}\xi^{k+s}}$ then  we get the following recurrence relation :\\
 
 \n $\displaystyle{ k(k-2)a_{k}  + i[\rho (k+1)+ \frac{ \sigma}{\lambda} ] a_{k+1} - 2(k+2)a_{k+2} = 0}; \, \, k \geq 2$ $\hfill { } {\color{blue}(3.20)}$\\
 
\n {\color{red}$\bullet$}  It is difficult to know if the sequence $(a_{k})_{k \in \mathbb{N}}$ belongs to $\mathbb{B}_{s}$.\\

\n {\color{red}$\bullet$}  In the vicinity of $ z = i \sqrt{2}\xi = \infty$ , the difficulty remains even if we use the techniques of the Thomé series in equation {\color{blue}(3.3)} by using the proposition $8$ of {\color{blue}[22]} where  one can choose as two linearly independent solutions of {\color{blue}(3.3)} the functions $\psi_{1}$ and $\psi_{2}$, defined by the power series :\\

\n $\displaystyle{\psi_{1}(\xi) = e^{-\frac{1}{2}(\gamma - \alpha - 2)} \sum_{n=0}^{\infty}a_{n}(\alpha, \beta, \gamma, \delta)\xi^{-n}}$  \hfill { } {\color{blue}(3.21)}\\

\n $\displaystyle{\psi_{2}(\xi) = e^{-\frac{1}{2}(\gamma + \alpha + 2}e^{\xi^{2} + \beta \xi} \sum_{n=0}^{\infty}a_{n}(\alpha, i\beta, -\gamma, - i\delta)(- i\xi)^{-n}}$  \hfill { } {\color{blue}(3.22)}\\

\n where $\displaystyle{\alpha = -1, \beta = - \frac{\mu \sqrt{2}}{\lambda}, \gamma = 1, \delta = \frac{2\sigma \sqrt{2}}{\lambda}}$  and the complex numbers $a_{n }$ are defined by {\color{blue}($\star_{3}$)}.\\

\n Now as  $\varphi(z) = \psi(\xi)$ where $z = i\sqrt \xi$,  then we deduce from {\color{blue}(3.21)} and {\color{blue}(3.22)} that {\color{blue}(3.2)}  at $\infty$, satisfies :\\

\n $\displaystyle{\varphi_{1}(z) \sim const.}$ as $z \longrightarrow \infty$  \hfill { } {\color{blue}(3.23)}\\

\n $\displaystyle{\varphi_{2}(z) \sim \frac{1}{z}e^{-\frac{1}{2}z^{2} + i\rho z}}$ as $z \longrightarrow \infty$  \hfill { } {\color{blue}(3.24)}\\

\n So that for any given direction in the complex $z-$plane an eigenfunction will be asymptotically a linear combination of {\color{blue}(3.23)} and {\color{blue}(3.24)}\\

\n In the following, we will show that  Eq {\color{blue}(3.24)} is not acceptable in Bargmann space and thus it can be rejected.\\

\n {\color{red}$ \bullet$} As  there are many important applications of the closedness of the range  of linear operator in the spectral study of differential operators and also in the context of perturbation theory (see e.g. {\color{blue}[Goldberg]}) {\color{blue}[9]} then to require an asymptotic behavior of the solutions of the eigenvalue problem associated with $H_{\mu, \lambda}$, we begin by showing that the image of $H_{\mu, \lambda}$ is closed that permit us to deduce a natural asymptotic behavior of the solutions of the eigenvalue problem associated to $H_{\mu, \lambda}$ in vicinity of $z = \infty$\\

\n {\bf{\color{red}Proposition 3.3}}\\

\n Let $\displaystyle{\mathbb{B}_{0} = \{ \varphi \in \mathbb{B} ; \varphi(0) = 0\} }$ and $H_{\mu, \lambda}^{mim} = H_{\mu, \lambda}$ with domain $D(H_{\mu, \lambda}^{mim})$ and $\mathcal{R}(H_{\mu, \lambda}^{mim})$  denotes the range of $H_{\mu, \lambda}^{mim}$ then \\

\n i) In  Bargmann space $\mathbb{B}$,  we have $\displaystyle{\varphi \in \mathbb{B} \iff  z \longrightarrow g(z) = \frac{\varphi^{'}(z) - \varphi^{'}(0)}{z} \in  \mathbb{B}}$.\\

\n ii) For $\mu \neq 0$, $\displaystyle{\mid\mid H_{\mu, \lambda}^{mim} \varphi  \mid\mid  \geq \mid \mu \mid \mid\mid \varphi \mid\mid , \, \forall \, \varphi \in D(H_{\mu, \lambda}^{mim})}$.\\

\n iii) $\mathcal{R}(H_{\mu, \lambda}^{mim})$ is closed in Bargmann space. \\

\n iv) $\mathcal{R}(H_{\mu, \lambda}^{mim})$ is dense in Bargmann space. \\

\n v) For $\mu \neq 0$,  $\displaystyle{H_{\mu, \lambda}^{mim}}$ is invertible. \hfill { } {\color{blue} $\square$}\\

\n {\bf {\color{red}Remark 3.4}}\\

\n 1) The properties of this proposition are well known. In particular (i) is the lemma $1$ page $127$ in {\color{blue}[3]} see also {\color{blue}[12]}  or lemma 0.3 page $338$ in {\color{blue}[15]}.\\

\n 2) It is convenient to reproduce the proofs of these properties which serve as a starting point of our study on the completeness of eigenfunctions of $H_{\mu, \lambda}$.\hfill { } {\color{blue} $\square$}\\

\n {\bf{\color{red}Proof}} \\

\n i) For $\displaystyle{\phi (z) = \sum_{n=0}^{\infty}a_{n}z^{n}}$, we have $\displaystyle{\phi (z) = \sum_{n=0}^{\infty}a_{n+1}z^{n+1}}$ and $\displaystyle{\frac{\phi' (z) - \varphi'(0)}{z} = \sum_{n=0}^{\infty}a_{n+2}z^{n+2}}$.\\

\n By using the isometry between spaces $\mathbb{B}$ and  $\mathbb{B}_{s}$, we deduce that \\

\n $\displaystyle{\mid\mid \frac{\phi' (z) - \varphi'(0)}{z} \mid\mid^{2} = \sum_{n=0}^{\infty}n!(n + 2)^{2} \mid a_{n+2}\mid^{2}}$.\\

\n i.e.\\

\n $\displaystyle{\mid\mid \frac{\phi' (z) - \varphi'(0)}{z} \mid\mid^{2} = \sum_{n=0}^{\infty}\frac{n+ 2}{n + 1}(n + 2)! \mid a_{n+2}\mid^{2}}$. $\hfill { } {\color{blue}(3.23)}$\\

\n As $\displaystyle{\frac{n + 2}{n + 1} \leq 2}$ then  $\displaystyle{\mid\mid \frac{\phi' (z) - \varphi'(0)}{z} \mid\mid^{2} \leq 2 \sum_{n=0}^{\infty}(n + 2)! \mid a_{n+2}\mid^{2}}$  Then \\

\n if $\varphi \in \mathbb{B}$, we have  $\displaystyle{ z \longrightarrow \frac{\varphi^{'}(z) - \varphi^{'}(0)}{z} \in  \mathbb{B}}$.\\

\n For the reciprocity, we remark that $\displaystyle{\frac{n + 2}{n + 1} > 1}$ and from   {\color{blue}(3.23)} we deduce that $\displaystyle{ \mid\mid \varphi \mid\mid \leq \mid\mid g \mid\mid}$ where $\displaystyle{g(z) = \frac{\varphi^{'}(z) - \varphi^{'}(0)}{z}}$.\\

\n ii) For $\lambda \in \mathbb{R}$  we observe that \\

\n $ \displaystyle{ \mid \Re e < H_{\mu, \lambda}^{min}\varphi , \varphi > \mid  = \mid \mu \mid\mid\mid A \varphi \mid\mid^{2} =  \mid \mu \mid \sum_{n=0}^{\infty}(n+1)\mid a_{n}\mid^{2}  \geq \mid\mu\mid \sum_{n=0}^{\infty}\mid a_{n}\mid^{2} = \mid \mu \mid\mid\mid \varphi \mid\mid^{2}}$ \hfill { } {\color{blue}(3.24)}\\

\n Now by using the Cauchy's inequality, we deduce that \\

\n {\color{red}($\alpha$)} $\displaystyle{\mid\mid H_{\mu, \lambda}^{min} \varphi  \mid\mid  \geq \mid \mu \mid\mid\mid A \varphi \mid\mid , \, \forall \, \varphi \in D(H_{\mu, \lambda}^{min} )}$  \hfill { } {\color{blue}(3.25)}\\

\n In particular\\

\n {\color{red}($\beta$)}  the injection of $\displaystyle{D(H_{\mu, \lambda}^{min})}$ in $D(A)$ is continuous. \hfill { } {\color{blue}(3.26)}\\

\n and\\

\n {\color{red}($\gamma$) } $\displaystyle{\mid\mid H_{\mu, \lambda}^{min} \varphi \mid\mid  \geq \mid \mu \mid \mid\mid \varphi \mid\mid , \, \forall \, \varphi \in D(H_{\mu, \lambda}^{min} )}$  \hfill { } {\color{blue}(3.27)}\\

\n (iii) From  {\color{blue}(3.25)} and as the injection of $D(A)$ in $\mathbb{B}$ is a compact mapping, we deduce that :\\

\n The injection $D(H_{\mu, \lambda}^{min})$  $\longrightarrow$ $\mathbb{B}$ is a compact mapping.  \hfill { } {\color{blue}(3.28)}\\

\n (iv) From  {\color{blue}(3.27)}, it follows that $\displaystyle{ \mathcal{N}( H_{\mu, \lambda}^{min}) = \{0\}}$ and  $\displaystyle{ \mathcal{R}( H_{\mu, \lambda}^{min})}$ is closed, where $\displaystyle{ \mathcal{N}( H_{\mu, \lambda}^{min})}$ and $\displaystyle{ \mathcal{R}( H_{\mu, \lambda}^{min})}$ denote the kernel and the range of $H_{\mu, \lambda}^{min}$ respectevly.\\

\n (v) In fact , $\displaystyle{< H_{\mu, \lambda}^{min} \varphi , \psi > = 0 \iff  - i\lambda z \psi^{''}(z)  + (- i \lambda z^{2} + \mu z)\psi^{'}(z) = 0  \iff H_{\mu, - \lambda}\psi = 0}$.\\

\n for $\lambda \neq 0$, we have $\displaystyle{\psi'(z) = e^{-\frac{1}{2} z^{2} + i \frac{\mu}{\lambda}z}}$.\\

\n Now we show that the normalizability requirement for $\psi(z)$ is not verified in some direction of the $z-$plane.\\

\n In fact, we will use the fundamental property (i) for $z \longrightarrow - i \infty$.\\

\n Let $z = \Re e z + i \Im m z$, then  $\displaystyle{e^{- \mid z \mid^{2}}\mid \frac{\psi'(z)}{z} \mid = e^{- \mid z \mid^{2}}\mid \frac{ e^{-\frac{1}{2} z^{2} + i \frac{\mu}{\lambda}z}}{z} \mid^{2} = \frac{1}{\mid z \mid^{2}} e^{-2(\Re e z)^{2}} e^{-2 \frac{\mu}{\lambda}\Im m z }}$.\\

\n Given $\mu$ and $ \lambda$, there exists a direction of $z-$plane  where  $\displaystyle{\psi(z) = \int_{0}^{z} e^{-\frac{1}{2} z\xi^{2} + i \frac{\mu}{\lambda}\xi}d\xi}$ cannot be considered in Bargmann space.\\

\n {\color{red}{\bf Remark 3.5}} \\

\n 1) To sum up  in $\displaystyle{\mathbb{B}_{0} = \{ \varphi \in \mathbb{B} ; \varphi(0) = 0 \}}$ the Hilbert space orthogonal to the vacuum (the only state of exactly zero energy), the eigenvalue conditions of $H_{\mu, \lambda}\varphi(z) = \sigma \varphi(z)$ are :\\

\n $\varphi(z) \sim z$  for $z \sim 0$  \hfill { } {\color{blue}(3.29$_{a}$) }\\

\n and \\

\n   $\varphi(z) \sim const.$  for $z \sim - i\infty$ \hfill { } {\color{blue}(3.29$_{b}$) }\\

\n  2) From iv) and v) we deduce that $H_{\mu, \lambda}^{min}$ is invertible i.e., \\

\n $0 \in \rho(H_{\mu, \lambda}^{min}$ where $\rho(H_{\mu, \lambda}^{min})$ denotes the resolvent set of $H_{\mu, \lambda}^{min}$.  \hfill { } {\color{blue}(3.30) }\\

\n 4) From {\color{blue}(3.28)} and as $\rho(H_{\mu, \lambda}^{min}) \neq \emptyset$ we deduce that the resolvent of $H_{\mu, \lambda}^{min}$ is compact. \\

\n 5) From {\color{blue}[15]}  (see theorem 1.2), it is well known that  for $\mu \neq 0$ we have $D(H_{\mu, \lambda}^{min}) = D(H_{\mu, \lambda})$.\\

\n  6) From  lemmas 9, 10 and 11  of  the ref. {\color{blue}[13]} (see page $280$), it is well known that the eigenvalues of $H_{\mu, \lambda}$ are real. \\

\n {\bf {\color{red} $\S$ 4 A transformation procedure of $H_{\mu, \lambda}$ to a symmetric operator with compact resolvent on the negative imaginary axis}}\\

\n We begin by recalling an original lemma on Bargmann space with its proof was established in {\color{blue}[13]}  \\

 \n {\bf{\color{red} Lemma 4.1}}  \\

\n Let $\mathbb{B}$ be Bargmann space,  if $\varphi \in \mathbb{B}$ then its restriction on  $x + i\mathbb{R}$ is square integrable function with measure $\displaystyle{e^{-\mid y \mid^{2}}dy , y \in \mathbb{R}}$ for all  fixed $x \in \mathbb{R}$.\\

\n {\bf Proof}\\

\n Bargmann has built an isometry between the space $\mathbb{B}$ and $L_{2}(\mathbb{R})$ so that,  $\forall \,\varphi \in \mathbb{B}$ is uniquely represented by $f \in  L_{2}(\mathbb{R})$ by means of the following integral:\\

\n $\displaystyle{\varphi (z) = c\int_{\mathbb{R}}e^{- \frac{z^{2}}{2} - \frac{q^{2}}{2} +\sqrt{2}zq } f(q)dq  }$  $\hfill { }  {\color{blue}(4.1)}$\\

\n such that $\displaystyle{\mid\mid \varphi \mid\mid = \mid\mid f \mid\mid_{L_{2}( \mathbb{R})}}$ $\hfill { }  {\color{blue}(4.2)}$\\

\n Now, we put $\displaystyle{g_{z}(q) = e^{- \frac{q^{2}}{2} +\sqrt{2}zq }}$, which is function belonging to $L_{2}(\mathbb{R})$. Then we can write {\color{blue}(4.1)} in the following form :\\

\n $\displaystyle{e^{\frac{z^{2}}{2}}\phi (z) = }c\displaystyle{\int_{\mathbb{R}}g_{z}(q)f(q)dq  }$ \hfill { } {\color{blue} (4.3)}\\

\n As $g_{z}(q)$ and $f$ are in  $L_{2}(\mathbb{R})$, we can apply the Parseval's identity at {\color{blue}(4.3)} to get :\\

\n $\displaystyle{e^{\frac{z^{2}}{2}}\varphi (z) = c_{1}}\displaystyle{\int_{\mathbb{R}}\hat{g}_{z}(p)\hat{f}(p)dp}$ \hfill { }  {\color{blue}(4.4)}\\

\n As $\displaystyle{g_{z}(q)}$ is a Gaussian's function,  one knows how to calculate its Fourier transform:\\

\n $\displaystyle{\hat{g}_{z}(p) = \int_{\mathbb{R}}e^{ipq}g_{z}(q)dq =\int_{\mathbb{R}}e^{- \frac{q^{2}}{2} + (\sqrt{2}z - ip)q} dq = (2\pi)^{\frac{1}{2}}e^{\frac{1}{2}(\sqrt{2}z - ip)^{2}}}$\\

 \n i.e.\\

 \n $\displaystyle{\hat{g}_{z}(p) = (2\pi)^{\frac{1}{2}}e^{\frac{1}{2}(\sqrt{2}z - ip)^{2}}}$\hfill { } {\color{blue} (4.5})\\

 \n Let $z = x + iy$, then for all fixed $x \in \mathbb{R}$, the function $\hat{g}_{z}(p)$ can be written under following form:\\

\n  $\displaystyle{\hat{g}_{z}(p) = (2\pi)^{\frac{1}{2}}e^{\frac{1}{2}(\sqrt{2}x - i(p - \sqrt{2}y)^{2}}}$ \hfill { }  {\color{blue}(4.6)}\\

\n  If we put $\displaystyle{h_{x}(p) = e^{x^{2}}e^{-\frac{1}{2}(p^{2} + i\frac{\sqrt{2}}{2}px)}}$ then $\displaystyle{e^{\frac{z^{2}}{2}}\varphi (z)}$ is the convolution product of $h_{x}$ with $\hat{f}$ evaluated in $\sqrt{2}y$, i.e.,\\

\n $\displaystyle{e^{\frac{(x + iy)^{2}}{2}}\varphi (x + iy) = C_{2} h_{x}*\hat{f}(\sqrt{2}y) }$\hfill { }  {\color{blue}(4.7)}\\

\n It follows that $\displaystyle{y \rightarrow e^{\frac{(x + iy)^{2}}{2}}\varphi (x + iy)}$ is in  $L_{2}(\mathbb{R})$ for all $x \in \mathbb{R}$\\.

\n By applying the Young's inequality we deduce that :\\

\n $\mid\mid e^{\frac{(x + iy)^{2}}{2}}\varphi (x + iy)\mid\mid^{2}_{L_{2}(\mathbb{R})} \leq \mid\mid h_{x}\mid\mid^{2}_{L_{1}(I\!\!R^{N})}.\mid\mid \hat{f}\mid\mid^{2}_{L_{2}(\mathbb{R})} \hfill { }  {\color{blue}(4.8)}$\\

\n where $C_{1}, C_{2}$ and $C$ are constants.\\

\n As $\mid\mid h_{x}\mid\mid^{2}_{L_{1}(\mathbb{R})} = C_{3}e^{2\mid x \mid^{2}}$ et $\mid\mid \hat{f}\mid\mid^{2}_{L_{2}(\mathbb{R})} = C_{4}\mid\mid f \mid\mid^{2}_{L_{2}(\mathbb{R})}$, then we deduce that :\\

\n $\displaystyle{\int_{\mathbb{R}}e^{-\mid y \mid^{2}}\mid \varphi (x + iy)\mid^{2}dy \leq C_{5}e^{-\mid x \mid^{2}}\mid\mid \varphi \mid\mid^{2}}$ $\hfill { }  {\color{blue}(4.9)}$ \\

\n where $C_{5}$ depends of  all the previous constants.\\

\n {\color{red}$\bullet$} From the above lemma, it is possible  to restrict ourselves to the negative imaginary axis, $z = - iy , y > 0$, where we impose the boundary conditions {\color{blue}(3.29)}.\\

\n Let us write therefore $\displaystyle{u(y) = \varphi(-iz)}$ then we have\\

\n  $\displaystyle{H_{_{R}}^{\mu, \lambda}u(y) = - \lambda y u''(y) + (\lambda y^{2} + \mu y)u'(y) = \lambda [-y\frac{d^{2}}{dy^{2}} + (y^{2} + \rho y)\frac{d}{dy}]u(y)}$ ; $\hfill { }  {\color{blue}(4.10)}$ \\

\n $\displaystyle{ \rho = \frac{\mu}{\lambda} }$\\

\n The eigenfunction equation $H_{_{R}}{\mu, \lambda}u(y) = \sigma u(y)$ may be rewritten (for $\lambda \neq 0$) in the form :\\

\n $\left \{\begin{array}{c} \displaystyle{\lambda y[-\frac{d^{2}}{dy^{2}} + (y + \rho )\frac{d}{dy}]u(y) = \sigma u(y)} \quad \quad\\
\quad\\
u(y )\sim y , \quad y \longrightarrow 0 \quad \quad \quad \quad \quad \quad \quad \quad \quad \\
\quad\\
u(y) \sim const. , \quad y \longrightarrow +\infty \quad \quad \quad \quad \quad  \quad \\
\end{array}\right.$ \hfill { }  {\color{blue}(4.11)}\\
\quad\\

\n The term proportional to $u'(y)$ can be eliminated by the transformation \\

\n $\displaystyle{u(y) =T(y) v(y)\,\, \mbox{where} \,\, T(y) = e^{\frac{1}{4}(y + \rho )^{2}}}$ \hfill { }  {\color{blue}(4.11)}\\

\n that gives the eigenvalue equation \\

\n $\left \{\begin{array}{c} \displaystyle{\lambda y[-\frac{d^{2}}{dy^{2}} + \frac{1}{4}(y + \rho )^{2} - \frac{1}{2}]v(y) =\sigma v(y)} \\
\quad\\
v(y) \sim y , \quad y \longrightarrow 0 \quad \quad \quad \quad \quad \quad \quad \quad \quad \quad\\
\quad\\
\displaystyle{v(y) \sim e^{-\frac{y}{2}( \frac{y}{2} + \rho)}}, \quad y \longrightarrow +\infty \quad \quad \quad \quad \quad  \quad \\
\end{array}\right.$ \hfill { }  {\color{blue}(4.12)}\\

\n As  in vicinity of $y = 0$  we have $v(y) \sim y$  and in vicinity of $y = +\infty$ we have $\displaystyle{\mid v(y) \mid^{2} \sim e^{-(\frac{y^{2}}{2} + \rho y)}}$ then  $v$ vanishes at origin and $v \in L^{2}[0, \infty[$.\\

\n We denote by $L_{R}^{2}[0, \infty[$ the space of functions in $L^{2}[0, \infty[$ satisfying both contraints.\\

\n Let  $\displaystyle{\underline{H}^{\mu, \lambda} = e^{-T(y)}H_{_{R}}^{\mu, \lambda}e^{T(y)} }$, then $\displaystyle{\underline{H}{\mu, \lambda} =  \lambda y[-\frac{d^{2}}{dy^{2}} + \frac{1}{4}(y + \rho )^{2} - \frac{1}{2}]}$.\\

\n  Now, if we consider  the change of variable $y = x^{2}$  and the similarity transform given by $y^{\frac{1}{4}} = x^{\frac{1}{2}}$ then we  get \\

\n $\displaystyle{\tilde{H}_{\mu, \lambda} = x^{-\frac{1}{2}}\underline{H}^{\mu, \lambda}x^{\frac{1}{2}}}$ \\

\n $\displaystyle{ = x^{-\frac{1}{2}} e^{-T(x^{2})} H_{_{R}}{\mu, \lambda} e^{T(x^{2})}x^{\frac{1}{2}}}$ \\

\n $\displaystyle{ = \frac{\lambda}{4} [-\frac{d^{2}}{dx^{2}} + \frac{3/4}{x^{2}} + x^{2}((x^{2} + \rho)^{2} - 2)] }$ \hfill { }  {\color{blue}(4.13)}\\

\n $\displaystyle{v(y) = (x^{2})^{\frac{1}{4}} w(x) =  x^{\frac{1}{2}} w(x) ,  \quad y = x^{2}}$.\\

\n In particular \\

\n $\displaystyle{w(x) = x^{-\frac{1}{2}} u(y) \sim x^{\frac{3}{2}} , \, \, x \longrightarrow 0}$  \hfill { }  {\color{blue}(4.14)}\\

\n  and \\

\n  $\displaystyle{w(x) \sim  x^{-\frac{1}{2}}e^{-\frac{1}{4}(x^{2} + \rho)^{2}}, \,\,  x \longrightarrow \infty}$.  \hfill { }  {\color{blue}(4.15)}\\

\n {\color{red}{\bf Remark 4.2}}\\

\n 1) The similarity transformation from $H_{\mu, \lambda}$ to $\tilde{H}_{\mu, \lambda}$ is {\color{red}non-unitary}, but {\color{blue}bijective}.\\

\n 2) On $L_{2}[0, +\infty[$,  $\tilde{H}_{\mu, \lambda}$  with domain $D(\tilde{H}_{\mu, \lambda}) = \{w  \in L_{2}[0, +\infty[; \tilde{H}_{\mu, \lambda}w \in L_{2}[0, +\infty[\} $ is \\

\n a symmetric operator.\\

\n 3)  As the injection of maximal domain of $ H_{\mu, \lambda}$ in Bargmann space $\mathbb{B}$ is compact then we deduce \\

\n  that:\\

\n (i) the injection of maximal domain of $\tilde{H}_{R}^{\mu, \lambda}$ in  \\

\n $\displaystyle{\mathbb{B}_{R} = \{u : [0, + \infty[ \longrightarrow \mathbb{C}\, \mbox{analytic}; \, \int_{0}^{\infty} \mid u(y)\mid^{2}e^{-y^{2}}dy < +\infty\}}$ \\

\n is compact.\\

\n  and \\

\n (ii) the injection of maximal domain of $\tilde{H}_{\mu, \lambda}$ in  $L_{2}[0, +\infty[$ is also compact, in particular $\tilde{H}_{\mu, \lambda}$ is a symmetric operator with compact resolvent.\\

\n 4) The equation $H_{\mu, \lambda} u = \sigma u$ is  equivalent to $\tilde{H}w = \sigma w$, where \\

\n $\displaystyle{w(x) = x^{-\frac{1}{2}}e^{-\frac{1}{4}(x^{2} + \rho)^{2}}u(y = x^{2})}$. \hfill { }  {\color{blue}(4.14)}\\

\n is an analytical function in $x$.\\

\n We end this section by the main result of this paper\\

\n {\color{red} {\bf Proposition 4.3}}\\

\n The eigenfunctions of $H_{\mu, \lambda} , \mu > 0$  form a complete basis in Bargmann space.\\

\n {\bf {\color{red}Proof }}\\

\n The above formulated problem of finding the spectrum of the Hamiltonian $H_{\mu, \lambda} $ or of $H_{_{R}}^{\mu, \lambda} $  acting on Bargmann space or on its  restriction functions respecti to negative imaginary axis is equivalent to the problem of finding the spectrum of the symmetric  operator $\tilde{H}_{\mu, \lambda}$
 in the space $L^{2}([0,+\infty[)$.\\

\n  As  $\tilde{H}_{\mu, \lambda}$ is symmetric with compact resolvent, then by using the spectral theorem for symmetric operators with compact resolvent, we deduce that the eigenfunctions of  $\tilde{H}_{\mu, \lambda}$ form a complete basis :\\

\n {\color{red}{\bf Theorem 4.4}} ({\color{blue} spectral theorem for symmetric operators with compact resolvent})\\

\n  1) Let $T$ be a densely defined linear operator  in a Hilbert space  $\mathcal{H}$ and let $\sigma \in \rho(T) \bigcap \mathbb{R}$ then $T$ is symmetric  iff $\displaystyle{R_{\sigma} = (T - \sigma I)^{-1}}$ is self-adjoint.\\

\n  2) Let $T$ be a symmetric operator in a Hilbert space $\mathcal{H}$ such that $\displaystyle{(T - \sigma_{0} I)^{-1}}$ is compact, then  
$\displaystyle{\sigma_{p}(T) = \{ \sigma_{0} + \frac{1}{\gamma}; \gamma \in \sigma_{p}(T - \sigma_{0} I)^{-1})\}}$ where $\sigma_{p}(T)$ denotes the ponctual spectrum of $T$. \\

\n  3) If $T$ is a symmetric operator with compact resolvent then there exists $\sigma_{1} \in  \rho(T) \bigcap \mathbb{R}$  such that \\

\n $\displaystyle{ (T - \sigma_{1} I)^{-1}}$ is a compact self-adjoint operator.\\

\n 4) Let $T$ be a symmetric operator with  compact resolvent in a Hilbert space $\mathcal{H}$.\\

\n  For each $\sigma \in  \sigma_{p}(T)$, define $n(\sigma) = dim \,ker(T  - \sigma I )$, and choose an orthonormal basis \\

\n $\displaystyle{\{\varphi_{\sigma,1}, \varphi_{\sigma,2}, ......., \varphi_{\sigma, n(\sigma)}\}}$ of $ker(T  - \sigma I )$. Then $\displaystyle{\{\varphi_{\sigma,i}; \sigma \in  \sigma_{p}(T), i=1, 2, ....., n(\sigma)\}}$ forms \\

\n a complete orthonormal set for $\mathcal{H}$.: Moreover, we have the following:  \\

\n  i) $\displaystyle{ \varphi  \in  D(T)  \iff  \sum_{\sigma  \in  \sigma_{p}(T)}  \sum_{i=1}^{n(\sigma)}  \mid  \sigma  \mid^{2}  \mid <  \varphi , \, \varphi_{\sigma, i} > \mid^{2} <  \infty}$.\\

\n  ii) For any $\displaystyle{ \varphi  \in  D(T),  T\varphi  =  \sum_{\sigma  \in  \sigma_{p}(T)}  \sum_{i=1}^{n(\sigma)}   \sigma <  \varphi , \, \varphi_{\sigma, i} > \varphi_{\sigma, i}}$. \\

\n {\bf{\color{red} $\S$ 5 \, An explicit integral form of the inverse of  $H_{\mu, \lambda}$ on negative  imaginary axis and some spectral properties of $H_{\mu, \lambda}^{-1}$}}\\

\n {\bf {\color{blue} A) Implementation of integral equation associated to $H_{\mu, \lambda}$}}\\

\n On Bargmann space, we consider the  equation $H_{\mu, \lambda} \varphi(z)  = \tilde{\sigma}  \varphi(z) $. If we put $\displaystyle{\rho = \frac{\mu}{\lambda}}$ and\\ $\displaystyle{\sigma = \frac{\tilde{\sigma}}{\lambda} }$ with $\lambda \neq 0$ then this eigenvalue equation can be written in the following form :\\

 \n $\displaystyle{\varphi^{"}(z) + p(z) \varphi^{'}(z) + q(z) \varphi(z) = 0}; \displaystyle{ p(z) = (z - i\rho)}$ and $\displaystyle{ q(z) =  \frac{i\sigma}{\lambda z}}$ $\hfill { } {\color{blue}(5.1)} $\\

\n We restrict the study of the equation {\color{blue}(5.1)} on the lines parallel to the imaginary axis :\\

\n Let $z = a + iy$ and $\varphi(a + iy) = \varphi_{a}(y)$ where $a$ is fixed. Then we get\\

\n $\left \{ \begin{array} [c] {l} \displaystyle{\varphi_{a}^{"}(y) - (y - \gamma) \varphi_{a}^{'}(y) - \frac{i\sigma}{a + iy}\varphi_{a}(y) = 0}\\
\quad\\
\gamma = \rho + ia\\
\end{array} \right . $ $\hfill { } {\color{blue}(5.2)}$
\quad\\

\n We put \\

\n $\displaystyle{\varphi_{a}^{'}(y) = K_{a}(y)e^{\frac{1}{2}(y - \gamma)^{2}}}$ $\hfill { } {\color{blue}(5.3)} $\\

\n and substitute it in {\color{blue}(5.2)} then we get\\

\n $\displaystyle{K_{a}^{'}(y) = i\sigma e^{-\frac{1}{2}(y - \gamma)^{2}}\frac{\varphi_{a}(y)}{a + iy}}$ $\hfill { } {\color{blue}(5.4)} $\\

\n We choose a primitive representation of the function $\displaystyle{ y \longrightarrow K_{a}^{'}(y)}$ as follows\\

\n $\displaystyle{K_{a}(y) = i\sigma \int_{-\infty}^{y}e^{-\frac{1}{2}(u - \gamma)^{2}}\frac{\varphi_{a}(u)}{a + iu}du + c_{a} \, ; y \in ]-\infty, 0]}$ $\hfill { } {\color{blue}(5.5)} $\\

\n where $c_{a}$ is a constant.\\

\n Consequently, equation {\color{blue}(5.3)} can be written as follows\\

\n $\displaystyle{\varphi_{a}^{'}(y) = i\sigma e^{\frac{1}{2}(y - \gamma)^{2}} \int_{-\infty}^{y}e^{-\frac{1}{2}(u - \gamma)^{2}}\frac{\varphi_{a}(u)}{a + iu}du  + c_{a}e^{\frac{1}{2}(y - \gamma)^{2}} }$ $\hfill { } {\color{blue}(5.6)} $\\

\n {\bf{\color{red} Lemma 5.1}}\\

\n For $\rho > 0$ the representation {\color{blue}(5.5)} has a sense.\\

\n {\bf{\color{red} Proof}}\\

\n For $\varphi \in \mathbb{B}$ we have $\displaystyle{\int_{\mathbb{R}^{2}}e^{-(x^{2} + y^{2})}\mid \varphi(x + iy) \mid^{2} dxdy < + \infty}$ then we apply the Fubini theorem to deduce that  \\

\n $\displaystyle{\int_{\mathbb{R}}e^{- y^{2}}\mid \varphi_{x}(y) \mid^{2} dy < + \infty}$ $\hfill { } {\color{blue}(5.7)} $\\

\n Now, as  $\displaystyle{ \int_{-\infty}^{y}e^{-\frac{1}{2}(u - \gamma)^{2}}\frac{\varphi_{a}(u)}{a + iu}du = e^{-\frac{1}{2}\gamma^{2}}\int_{-\infty}^{y}e^{-\frac{1}{2}u^{2}}\varphi_{a}(u)\frac{e^{\gamma u}}{a + iu}du }$ and we look for $\varphi$ in Bargmann space $\mathbb{B}$, we have on the one hand $\displaystyle{u \longrightarrow  e^{-\frac{1}{2}u^{2}}\varphi_{a}(u)}$ is square integrable and on the other hand for $\rho > 0$ the function $\displaystyle{u \longrightarrow  \frac{e^{\gamma u}}{ a + iu}}$ is also square integrable on $]-\infty, 0]$. Consequently, $\displaystyle{ \int_{-\infty}^{y}e^{-\frac{1}{2}(u - \gamma)^{2}}\frac{\varphi_{a}(u)}{a + iu}du}$ has a sense. \\

\n {\bf{\color{red} Remark 5.2}}\\

\n (i) If $\rho < 0$, we choose a primitive representation of the function $\displaystyle{ y \longrightarrow K_{a}^{'}(y)}$ as follows\\

\n $\displaystyle{K_{a}(y) = i\sigma \int_{y}^{+\infty}e^{-\frac{1}{2}(u - \gamma)^{2}}\frac{\varphi_{a}(u)}{a + iu}du + c_{a} }$, $c_{a}$ is constant $\hfill { } {\color{blue}(5.8)}$\\

\n This representation  has  a sense.\\

\n (ii) From {\color{blue}(5.3)} and {\color{blue}(5.5)} we deduce that\\

\n $\displaystyle{\varphi_{a}^{'}(y) = i\sigma e^{\frac{1}{2}(y - \gamma)^{2}} \int_{y}^{+\infty}e^{-\frac{1}{2}(u - \gamma)^{2}}\frac{\varphi_{a}(u)}{a + iu}du + c_{a}e^{\frac{1}{2}(u - \gamma)^{2}}}$ $\hfill { } {\color{blue}(5.9)}$\\

\n {\bf{\color{red}Theorem 5.3}}\\

\n (i) In the representation {\color{blue}(5.5)} we have $c_{a} = 0$\\

\n (ii)  \n $\displaystyle{Lim \, e^{-\frac{1}{2}(y - \gamma)^{2}}\varphi_{a}^{'}(y) = 0}$ as $y \longrightarrow -\infty$ (boundary condition) $\hfill { } {\color{blue}(5.10)} $\\

\n {\bf{\color{red}Proof}}\\

\n (i) By multiplying the two members of the equation {\color{blue}(5.9)} by $\displaystyle{e^{-\frac{1}{2}(u - \gamma)^{2}}}$, we get\\

\n $\displaystyle{e^{-\frac{1}{2}(u - \gamma)^{2}}\varphi_{a}^{'}(y) = i\sigma  \int_{y}^{+\infty}e^{-\frac{1}{2}(u - \gamma)^{2}}\frac{\varphi_{a}(u)}{a + iu}du + c_{a}}$ $\hfill { } {\color{blue}(5.11)} $\\

\n As the function $\displaystyle{u \longrightarrow e^{-\frac{1}{2}(u - \gamma)^{2}}\frac{\varphi_{a}(u)}{a + iu}}$ is integrable we deduce that \\

\n $\displaystyle{Lim \, e^{-\frac{1}{2}(y - \gamma)^{2}}\varphi_{a}^{'}(y) = c_{a}}$ as $y \longrightarrow - \infty$ $\hfill { } {\color{blue}(5.12)} $\\

\n Now, if $c_{a} \neq 0$ then we get that \\

\n $\displaystyle{\varphi_{a}^{'}(y) \thicksim c_{a}e^{\frac{1}{2}(y - \gamma)^{2}}}$ and   $\displaystyle{\frac{e^{-\frac{1}{2}y^{2}}}{y}\varphi_{a}^{'}(y) \thicksim c_{a}\frac{e^{- \rho y}}{y}e^{\frac{\gamma^{2}}{2} - ia}}$ $\hfill { } {\color{blue}(5.13)} $\\

\n We consider the space\\

\n $\displaystyle{\mathbb{B}_{a} = \{ \varphi_{a} : \mathbb{R} \longrightarrow \mathbb{C} \,\, \mbox{entire} \,; \int_{\mathbb{R}}e^{-\frac{1}{2}y^{2}}\mid \varphi_{a}(y)\mid^{2}dy < +\infty\}}$ $\hfill { } {\color{blue}(5.14)} $\\

\n Then \\

\n ($\alpha$) From lemma {\color{blue}(4.1)}, we deduce  that $\forall \, \varphi \in \mathbb{B}$ then $\varphi_{a} \in \mathbb{B}_{a}$\\

\n ($\beta$) From (i) of proposition {\color{blue}(3.3)}, we deduce  that $\forall \, \varphi_{a} \in \mathbb{B}_{a}$ then the function \\

\n $\displaystyle{y \longrightarrow  \frac{\varphi_{a}^{'}(y) - \varphi_{a}^{'}(0)}{y} \in \mathbb{B}_{a}}$

\n As in {\color{blue}(5.13)} the function $\displaystyle{y \longrightarrow \frac{e^{- \rho y}}{y}e^{\frac{\gamma^{2}}{2} - ia}}$ is not square integrable then $\varphi_{a} \notin \mathbb{B}_{a}$\\

\n which is contradictory with the property ($\beta$) and therefore $c_{a} = 0$.\\

\n (ii) we apply (i) to  {\color{blue}(5.12)}.\\

\n {\bf{\color{red} Remark 5.4}}\\

\n By applying the above theorem, equation  {\color{blue}(5.6)} can be written as follows\\

\n $\displaystyle{\varphi_{a}^{'}(y) = i\sigma e^{\frac{1}{2}(y - \gamma)^{2}} \int_{-\infty}^{y}e^{-\frac{1}{2}(u - \gamma)^{2}}\frac{\varphi_{a}(u)}{a + iu}du; \, y \in \,]-\infty, 0] }$ $\hfill { } {\color{blue}(5.15)} $\\

\n {\bf {\color{blue}  B) Justification for crossing to the limit in the equation {\color{red}(5.15)} as $a \longrightarrow 0$}}\\

\n {\bf{\color{red} Proposition 5.5}}\\

\n If $a \longrightarrow 0$ in {\color{blue}(5.15)} then we get \\

\n $\displaystyle{\varphi_{0}^{'}(y) = \sigma e^{\frac{1}{2}(y - \rho)^{2}} \int_{-\infty}^{y}e^{-\frac{1}{2}(u - \rho)^{2}}\frac{\varphi_{0}(u)}{u}du; \, y \in \,]-\infty, 0] }$ $\hfill { } {\color{blue}(5.16)} $\\

\n {\bf{\color{red} Proof}}\\

\n As the functions $\displaystyle{a \longrightarrow \varphi_{a}^{'}(y)}$ and $\displaystyle{a \longrightarrow \frac{e^{-\frac{1}{2}(u - \rho - ia)^{2}}}{a + iu}\varphi_{a}(u)}$ are continuous  then $\displaystyle{Lim\, \varphi_{a}^{'}(y) = \varphi_{0}^{'}(y)}$ and $\displaystyle{Lim \,  \frac{e^{-\frac{1}{2}(u - \rho - ia)^{2}}}{a + iu}\varphi_{a}(u) =  \frac{e^{-\frac{1}{2}(u - \rho)^{2}}}{iu}\varphi_{0}(u)}$ as $a \longrightarrow 0$. $\hfill { } {\color{blue}(5.17)}  $\\

\n We have to show that $\displaystyle{u \longrightarrow \frac{e^{-\frac{1}{2}(u - \rho)^{2}}}{iu}\varphi_{0}(u)}$ is integrable. For this purpose we consider the following function\\

\n $\displaystyle{f_{a}(u) = \frac{e^{-\frac{1}{2}(u - \rho - ia)^{2}}}{a + iu}\varphi_{a}(u)}$ $\hfill { } {\color{blue}(5.18)} $\\

 \n  As for $\epsilon > 0$ enough small  we have $\varphi_{a}(a + iu) \thicksim a +iu ; u \in ]-\epsilon, 0]$ (by using  lemma {\color{blue}4.1}), we deduce that \\
 
 \n $\displaystyle{\mid f_{a}(u)\mid \leq e^{\rho u}}$. $\hfill { } {\color{blue}(5.19)} $\\
  
\n Now, for $u \in ]-\infty , -\epsilon[$, we observe that \\

\n $\displaystyle{\mid f_{a}(u)\mid \leq e^{-\frac{\rho^{2}}{2}}\frac{e^{\rho u}}{u}}$. $\hfill { } {\color{blue}(5.20)}$\\

\n As the function\\

\n $g(u) = \left \{ \begin{array}[c] {l} \displaystyle{e^{-\frac{\rho^{2}}{2}}\frac{e^{\rho u}}{u}}; u \in ]-\infty , -\epsilon[\\
\quad\\
\displaystyle{e^{\rho u}}; u \in ]-\epsilon, 0]\\
\end{array} \right .$  $\hfill { } {\color{blue}(5.21)} $\\

\n is integrable and the the sequence $f_{a}(u)$ converges pointwise to  $\displaystyle{\frac{e^{-\frac{1}{2}(u - \rho^{2}}}{iu}\varphi_{0}(u)}$ as $a \longrightarrow 0$ then by applying the Lebesgue's Dominated Convergence Theorem, we deduce that\\

\n $\displaystyle{Lim\, \int_{-\infty}^{y}\frac{e^{-\frac{1}{2}(u - \rho - ia)^{2}}}{a + iu}\varphi_{a}(u)du = \int_{-\infty}^{y}\frac{e^{-\frac{1}{2}(u - \rho^{2}}}{iu}\varphi_{0}(u)}$ as $a \longrightarrow 0$ $\hfill { } {\color{blue}(5.22)} $\\

\n Consequently, we get\\

\n $\displaystyle{\varphi_{0}^{'}(y) = \sigma e^{\frac{1}{2}(y - \rho)^{2}} \int_{-\infty}^{y}e^{-\frac{1}{2}(u - \rho)^{2}}\frac{\varphi_{0}(u)}{u}du; \, y \in \,]-\infty, 0] }$ $\hfill { } {\color{blue}(5.23)}$\\ 

\n {\bf{\color{red} Theorem 5.6}}\\

\n (i) An operator integral associated to $H_{\mu, \lambda}$ is given as follows\\

\n $\displaystyle{K\psi(y) = \int_{-\infty}^{0}\mathcal{N}(y, s)\psi(s)ds}$$\hfill { } {\color{blue}(5.24)} $\\ 

\n where\\

\n ($\alpha$) $\displaystyle{\psi(y) = \frac{\varphi_{0}^{'}(y)}{y}e^{-\frac{y^{2}}{2}}\theta(y)}$$\hfill { } {\color{blue}(5.25)} $\\ 

\n with $\theta(y) = \left \{ \begin{array}[c] {l}\,  y ; \quad y \in [-1, 0]\\
\quad\\
- 1; \quad y \in ]-\infty, -1]\\
\end{array} \right .$ \\
\quad\\

\n ($\beta$) $\displaystyle{\mathcal{N}(y, s) = e^{-\rho y}\frac{\theta(y)}{y}.\frac{s}{\theta(s)}e^{\frac{s^{2}}{2}}\int_{-\infty}^{min(y, s)}e^{-\frac{1}{2}(u - \rho)^{2}}{u}du}$ $\hfill { } {\color{blue}(5.26)} $\\ 

\n where $min(y, s) = \left \{ \begin{array}[c] {l} y ; \quad y \leq s\\
\quad\\
s; \quad  s \leq y \\
\end{array} \right .$ \\
\quad\\

\n ($\gamma$) $\displaystyle{\psi \in L_{2}(]-\infty, 0], \theta(y)dy)}$ $\hfill { } {\color{blue}(5.27)}$\\ 

\n (ii) For $\psi \in \mathbb{B}_{0}$, an explicit inverse of $H_{\mu, \lambda}$ restricted on imaginary axis ; $y \in [0,  +\infty[$ is given by\\

\n $\displaystyle{H_{\mu, \lambda}^{-1} \psi(-iy) = \int_{0}^{\infty}\mathcal{N}_{\mu,\lambda}(y, s)\psi(-is)ds}$ $\hfill { } {\color{blue}(5.28}) $\\

\n where  $\displaystyle{\mathcal{N}_{\mu,\lambda}(y, s) = \frac{1}{\lambda s}e^{-\frac{s^{2}}{2} - \frac{\mu}{\lambda}s}\int_{0}^{min(y, s)}e^{\frac{u^{2}}{2} +\frac{\mu}{\lambda}u}du}$\\

\n {\bf{\color{red} Proof}}\\

\n (i) We substitue $\displaystyle{\varphi_{0}(u) = -\int_{u}^{0} \varphi_{0}^{'}(s)ds}$ in {\color{blue}(5.23)} to obtain\\

\n $\displaystyle{\varphi_{0}^{'}(y) = -\sigma e^{\frac{1}{2}(y - \rho)^{2}} \int_{-\infty}^{y}\frac{e^{-\frac{1}{2}(u - \rho)^{2}}}{u}\int_{u}^{0}\varphi_{0}^{'}(s)ds; \, y \in \,]-\infty, 0] }$ $\hfill { } {\color{blue}(5.29)} $\\ 

\n By applying the Fubini theorem we get\\

\n $\displaystyle{\int_{-\infty}^{y}\frac{e^{-\frac{1}{2}(u - \rho)^{2}}}{u}\int_{u}^{0}\varphi_{0}^{'}(s)ds = \int_{-\infty}^{0}[\int_{-\infty}^{min(y, s)}\frac{e^{-\frac{1}{2}(u - \rho)^{2}}}{u}du] \varphi_{0}^{'}(s)ds}$$\hfill { } {\color{blue}(5.30)} $\\ 

\n and\\

\n $\displaystyle{\varphi_{0}^{'}(y) = -\sigma e^{\frac{1}{2}(y - \rho)^{2}}\int_{-\infty}^{0}[\int_{-\infty}^{min(y, s)}}$ $\displaystyle{\frac{e^{-\frac{1}{2}(u - \rho)^2}}{u}du] \varphi_{0}^{'}(s)ds}$$\hfill { } {\color{blue}(5.31)} $\\ 

\n or\\

\n $\displaystyle{e^{-\frac{y^{2}}{2}}\varphi_{0}^{'}(y)}$ $\displaystyle{ = -\sigma e^{\frac{\rho^{2}}{2}}\int_{-\infty}^{0}[ e^{\rho y} \int_{-\infty}^{min(y, s)}\frac{e^{-\frac{1}{2}(u - \rho)^{2}}}{u}du] \varphi_{0}^{'}(s)ds}$$\hfill { } {\color{blue}(5.32)} $\\ 

\n Now, if we put  $\displaystyle{\psi(y) = \frac{\varphi_{0}^{'}(y)}{y}e^{-\frac{y^{2}}{2}}\theta(y)}$$\hfill { } {\color{blue}(5.33)} $\\  

\n with $\theta(y) = \left \{ \begin{array}[c] {l}\,  y ; \quad y \in [-1, 0]\\
\quad\\
- 1; \quad y \in ]-\infty, -1]\\
\end{array} \right .$ \\
\quad\\

\n and if we substitute $\psi$ in {\color{blue}(5.32)}) then we get \\

\n $\displaystyle{\psi(y) = -\sigma e^{\frac{\rho^{2}}{2}}\int_{-\infty}^{0}\mathcal{N}(y, s) \psi(s)ds }$ $\hfill { } {\color{blue}(5.34)} $\\

\n Where $\mathcal{N}(y, s)$ is given by  {\color{blue}(5.27)}.\\

\n (ii) Let $\varphi \in \mathbb{B}_{0}$ and  $\psi \in \mathbb{B}_{0}$, we consider the equation $H_{\mu, \lambda}\varphi = \psi (z)$ i.e.\\

\n $\displaystyle{i\lambda z \varphi^{"}(z)  + (i\lambda z^{2} + \mu z)\varphi'(z) = \psi(z)}$ $\hfill { } {\color{blue}(5.35)} $\\

\n Let $ z = -iy$, $u(y) = \varphi(-iy)$ and $ f(y) = \psi(-iy)$ with $y \in [0, \infty[$ then {\color{blue}(5.34)} can be \\

\n written in the following form:\\

\n $\displaystyle{ - \lambda y u''(y) + (\lambda y^{2} + \mu y) u'(y) = f(y)}$ $\hfill { } {\color{blue}(5.36)} $\\

\n or\\

\n $\displaystyle{ - u''(y) + (y + \rho ) u'(y) = \frac{f(y)}{\lambda y}}$ $\hfill { } {\color{blue}(5.37)} $\\

\n with help lemma {\color{blue}(4.1)}, the property i) of proposition {\color{blue}(3.3)} and a similar technique used in \\ proof of (i) of this theorem we get {\color{blue}(5.28)}.\\

\n In the following, we show that the integral operator $K$ is an operator of Hilbert-Schmidt\\ on $\displaystyle{L_{2}(]-\infty, 0], \theta(y)dy)}$ and we give some properties of kernel $\displaystyle{\mathcal{N}(y, s)}$ (see {\color{blue}(5.28)}) $ \mbox{in}$ \\
$\displaystyle{ L^{2}(]-\infty, 0]\times]-\infty, 0], \theta(y)\theta(s)dyds)}$.\\

\n {\bf {\color{red}Theorem 5.7}}\\

\n Let $ \mathcal{N}(y, s)$ = $ e^{-\rho y} \frac{\theta(y)}{y} \frac{s}{\theta(s)}e^{\frac{1}{2}s^{2}}$ $\displaystyle{\int_{-\infty}^{min(y,s)}}$ $\frac{e^{-\frac{1}{2}(u - \rho)^{2}}}{u}du$ be the kernel of the equation  :\\

\n $\displaystyle{ - \lambda y u''(y) + (\lambda y^{2} + \mu y) u'(y) = \psi(-iy)}$ , $\psi \in \mathbb{B}_{0}$ $\hfill { } {\color{blue}(5.38)} $\\

\n then:\\

\n (i) $ \mathcal{N}(y, s)$ is Hilbert-Shmidt on $\displaystyle{\mathcal{N}(y, s)}$ in $\displaystyle{L^{2}(]-\infty, 0]\times]-\infty, 0], \theta(y)\theta(s)dyds)}$ .\\

\n (ii) $\displaystyle{\mathcal{N}(y, s)}$ in $\displaystyle{L^{2}(]-\infty, 0]\times]-\infty, 0], \theta(y)\theta(s)dyds)}$.\\

\n {\bf {\color{red}Proof}}\\

\n (i) To apply the theorem of dominated convergence, we will increase $\mathcal{N}(y, s)$ by a function $\tilde{\mathcal{N}}(y, s) \in L^{2}(]-\infty, 0]\times]-\infty, 0], \theta(y)\theta(s)dyds)$\\

\n For this we pose $min(y, s) = m$  then we have :\\

\n $\displaystyle{\mid \frac{e^{-\frac{1}{2}(u - \rho)^{2}}}{u}\mid \leq e^{-\frac{1}{2}\rho^{2}}\frac{1}{m(m - \rho)}\mid (\rho - u) e^{-\frac{1}{2}u^{2} + \rho u}\mid}$\\

\n The function $ u \longrightarrow \displaystyle{\frac{1}{u(u - \rho)}} $ is bounded for $-\infty < u \leq m < 0$\\

\n Hence a first increase:\\

\n $\mid \mathcal{N}(y, s) \mid \leq \displaystyle{e^{-\frac{1}{2}\rho^{2}}\mid \frac{e^{-\rho y}}{y}\frac{\theta(y)}{\theta(s)}se^{\frac{1}{2}s^{2}}\frac{e^{-\frac{1}{2}m^{2} + \rho^{m}}}{m(m - \rho)}\mid}$\\

\n {\color{red}$\bullet$} If $y \leq s$, we have: \\

\n $ m = y$ et $\displaystyle{e^{\frac{1}{2}s^{2} - \frac{1}{2}y^{2}} < 1}.$\\

\n Then :\\

\n $\mid \mathcal{N}(y, s) \mid \leq \displaystyle{e^{-\frac{1}{2}\rho^{2}}\mid \frac{\theta(y)}{y}\frac{s}{\theta(s)}\frac{1}{y(y - \rho)}}\mid$\\

\n  {\color{red}$\bullet$} If $s \leq y$, we have: \\

\n $ m = s$ et $e^{\rho (s - y)} < 1.$\\

\n Therefore\\

\n $\mid \mathcal{N}(y, s) \mid \leq e^{-\frac{1}{2}\rho^{2}}\mid \frac{\theta(y)}{y}\frac{s}{\theta(s)}\frac{1}{s(s - \rho)}\mid$\\

\n Finally we put :\\

\n $\tilde{\mathcal{N}}(y, s) = \left\{
  \begin{array}{ c }
e^{-\frac{1}{2}\rho^{2}}\mid \frac{\theta(y)}{y}\frac{s}{\theta(s)}\frac{1}{y(y - \rho)} \quad pour \quad y \leq s \\
\quad\\
e^{-\frac{1}{2}\rho^{2}}\mid \frac{\theta(y)}{y}\frac{s}{\theta(s)}\frac{1}{s(s - \rho)} \quad pour \quad s \leq y \\
\end{array} \right .$\\

\n to get  $\mid \mathcal{N}(y, s) \mid \leq \mid \tilde{\mathcal{N}}(y, s) \mid$\\

\n It remains to show that $ \tilde{\mathcal{N}}(y, s) \in L^{2}(]-\infty, 0]\times]-\infty, 0], \theta(y)\theta(s)dyds)$.\\

\n Let $\Delta = ]-\infty, 0]\times]-\infty, 0] -]-1, 0]\times]-1, 0], \Delta_{1} = \{(y,s) \in \Delta ; y \leq s\}$ and \\

\n  $\Delta_{2} = \{(y,s) \in \Delta ; s \leq y\}$.\\

\n Then\\

\n $\displaystyle{\int_{-\infty}^{0}\theta(y)dy}\displaystyle{\int_{-\infty}^{0}\theta(s)\mid \tilde{\mathcal{N}}(y, s) \mid^{2} ds} = $\\

\n $\displaystyle{\int_{\Delta_{1}}\theta(y)\theta(s)\mid \frac{\theta(y)}{y}\frac{s}{\theta(s)}\frac{1}{y(y - \rho)}\mid^{2}dyds} + \displaystyle{\int_{\Delta_{2}}\theta(y)\theta(s)\mid \frac{\theta(y)}{y}\frac{s}{\theta(s)}\frac{1}{s(s - \rho)}\mid^{2}dyds}$\\

\n As $\rho > 0$  and the function  $ u \rightarrow \frac{1}{(u - \rho)^{2}}$  is integrable at the origin, we deduce that:\\

\n $\displaystyle{\int_{-\infty}^{0}\theta(y)dy}\displaystyle{\int_{-\infty}^{0}\theta(s)\mid \tilde{\mathcal{N}}(y, s) \mid^{2} ds < +\infty}$\\

\n then $\mathcal{N}(y, s)$ is a kernel of Hilbert-Shmidt.\\

\n {\bf {\color{red} Corollary 5.8}} ({\color{blue}some spectral properties of $K$})\\

\n I) We note that the functions :\\ 

\n $ y \longrightarrow \displaystyle{\int_{-\infty}^{y}\frac{e^{-\frac{1}{2}(u - \rho)^{2}}}{u}du} ;\quad y \leq s$\\

\n and\\

\n  $s \longrightarrow \displaystyle{\int_{-\infty}^{s}\frac{e^{-\frac{1}{2}(u - \rho)^{2}}}{u}du} ;\quad s \leq y$ .\\

\n  are decreasing and tend to zero respectively when $ y \longrightarrow -\infty $ and $ s \longrightarrow -\infty $.\\

\n ii) The kernel $-\mathcal{N}(y, s)$ is non negative.\\

\n iii) the operator $K\psi(y) = \displaystyle{\int_{-\infty}^{0} \mathcal{N}(y, s)\psi(s)ds}$ is an operator of Hilbert-Shmidt.\\

\n (iv) As $K$ is a Hilbert-Schmidt operator then its spectrum is discrete and as its kernel is non negative then by using Jentzsch theorem {\color{blue}[17]} or Krein-Rutman theorem {\color{blue}[18]}, we deduce the existence of real eigenvalue of $K$. \\

\n In particular, we deduce  the existence of real eigenvalue of $H_{\mu, \lambda} ; \mu > 0$ and that $H_{\mu, \lambda}^{-1}$  is  a positive operator on Bargmann space $\mathbb{B}$ in the following sense:\\

\n Let $\mathcal{C} \subset \mathbb{B}$ be a closed cone with non empty interior $ \displaystyle{int(\mathcal{C})}$ then $\displaystyle{H_{\mu, \lambda}^{-1}(\mathcal{C}-\{0\}) \subset  int(\mathcal{C})}$.\\

\n \n (v) But the question of the existence  of $\sigma > 0$ such that $\displaystyle{(H_{\mu, \lambda} - \sigma I) ^{-1}}$ is a positive operator is always an open question.\\

\n {\color{red}{\bf References}} \\

\n {\color{blue}[1]}  Abramovitz, M. and  Stegun, I. A. :. Handbook of Mathematical Functions, New York (1968).\\

\n {\color{blue}[2]}  Agarwal, R.P. and Regan, D.O. :. Ordinary and Partial Differential Equations: With Special Functions, Fourier Series, and Boundary Value Problems, Lecture 5.\\

\n {\color{blue}[3]}  Ando, T. and Zerner, M. :. Sur une valeur propre d'un opérateur, Comm. Math. Phys. 93 (1984). \\

\n {\color{blue}[4]}  Bargmann, V. :. On a Hilbert space of analytic functions and an associated integral transform I, Comm. Pure
Appl. Math. 14 (1962) 187-214. \\

\n{\color{blue}[5]}   Batola, F. :. Une généralisation d'une formule de Meixner-tricomi, Can. J. Math., Vol. XXXIV, No. 2, 1982, pp. 411- 422\\

\n {\color{blue}[6]}  Bender,C. M., Steven A. and Orszag, S. A.  :. Advanced Mathematical Methods for Scientists and Engineers (McGraw-Hill Book Company, New York, 1978).\\

\n {\color{blue}[7]}   Bateman, H. and Erdelyi, A. :. Higher Transcendental Functions, Vol. 1, New York (1953).\\

\n {\color{blue}[8]}  Decarreau, A.  Emamirad, H. and  Intissar, A. :. Chaoticité de l'opérateur de Gribov dans l'espace de Bargmann,
C. R. Acad. Sci. Paris Sér. I Math. 331 (2000) 751-756.\\

\n {\color{blue}[9]}  Goldberg, S. :. Unbounded Linear Operators. Mc Graw Hill, New York, 1966.\\

\n {\color{blue}[10]}  Hille,  E. :. Ordinary Differential Equations in the Complex Domain, Dover Publications Inc , New edition (1997) \\

\n  {\color{blue}[11]}  Ince,  E.L. :. Ordinary Differential Equations, Dover, New York, 1956. \\

\n  {\color{blue}[12]}  Intissar, A., Le Bellac, M. and Zerner, M. :. Properties of the Hamiltonian of Reggeon field theory, Phys. Lett. B 113 (1982) 487-489.\\

\n  {\color{blue}[13]}   Intissar, A. :. Etude spectrale d'une famille d'opérateurs non-symétriques intervenant dans la théorie des
champs de Reggeons, Comm. Math. Phys. 113 (1987) 263-297. \\

\n  {\color{blue}[14]}   Intissar, A. :. Analyse de Scattering d'un opérateur cubique de Heun dans l'espace de Bargmann, Comm.
Math. Phys. 199 (1998) 243-256.\\

\n  {\color{blue}[15]}  Intissar, A. :. Spectral analysis of non-selfadjoint Jacobi-Gribov operator and asymptotic analysis of its generalized eigenvectors, Advances in Mathematics (China), Vol.44, no 3, (2015), 335-353, doi: 10.11845/sxjz.2013117b\\

\n  {\color{blue}[16]}  Intissar, A. and Intissar, J. K. :. On Chaoticity of the Sum of Chaotic Shifts with Their Adjoints in Hilbert Space and Applications to Some Weighted Shifts Acting on Some Fock-Bargmann Spaces, Complex Anal. Oper. Theory, Volume 11, issue 3, (2017), 491-505 \\

\n  {\color{blue}[17]}   Jentzsch, R :. Uber Integralgleichungen mit positiven Kern, J. Reine Angew. Math. 141 (1912), 235-244.\\

\n  {\color{blue}[18]} Krein, M. G. and M. A. Rutman, M. A. :. Linear operators leaving invariant a cone in a Banach space, Amer. Math. Soc. Transl. Ser. 1 10 (1950), 199-325 [originally Uspekhi Mat. Nauk 3 (1948), 3-95]\\

\n  {\color{blue}[19]}  Kristensson, G. :. Second Order Differential Equations: Special Functions and Their Classification (Springer Science, New York, 2010).\\

\n  {\color{blue}[20]}  Maroni, P. :. Biconfluent Heun equation, In: A. Ronveau (Ed.), Heun's diff erential equations, 191-249, Oxford University Press, Oxford 1995\\

\n  {\color{blue}[21]}  Morse, P.M. and H Feshbach, H. :. Methods of Theoretical Physics, I, {\bf{\color{blue} $http://people.physics.tamu.edu/pope/methods_pages2011$.pdf}}\\

 \n  {\color{blue}[22]}  Roseau, A. :. On the solutions of the biconfluent Heun equations, Bull. Belg. Math. Soc. Simon Stevin ,Volume 9, Number 3 (2002), 321-342. doi:10.36045/bbms/1102715058. https://projecteuclid.org/euclid.bbms/1102715058\\

\n   {\color{blue}[23]}  Yosida, K. :. Lectures on Differential and Integral Equations, pp. 37-40.\\
\end{document}